\documentclass[12pt]{amsart}
\usepackage{amscd,amssymb}
\usepackage{tikz}
\usepackage{tabularx} 
\usepackage{tabu}
\usepackage[plainpages,backref,urlcolor=blue]{hyperref}
\usepackage[all]{xy}

\theoremstyle{plain}
\newtheorem{thm}[subsection]{Theorem}
\newtheorem{lem}[subsection]{Lemma}
\newtheorem{prop}[subsection]{Proposition}
\newtheorem{cor}[subsection]{Corollary}

\theoremstyle{definition}
\newtheorem{rk}[subsection]{Remark}
\newtheorem{definition}[subsection]{Definition}
\newtheorem{ex}[subsection]{Example}
\newtheorem{conj}[subsection]{Conjecture}

\numberwithin{equation}{section}
\newcommand{\Pc}{{\mathcal P}}

\newcommand{\I}{{\mathcal I}}

\newcommand{\A}{{\mathcal A}}
\newcommand{\B}{{\mathcal B}}

\newcommand{\CC}{{\mathcal C}}
\newcommand{\LL}{{\mathcal L}}

\newcommand{\al}{{\alpha}}
\newcommand{\be}{{\beta}}

\newcommand{\SSS}{{\mathcal S}}

\newcommand{\Q}{\mathbb{Q}}

\newcommand{\C}{\mathbb{C}}
\newcommand{\PP}{\mathbb{P}}

\DeclareMathOperator{\mult}{mult}

\DeclareMathOperator{\ord}{ord}
\DeclareMathOperator{\mdr}{mdr}

\begin{document}

\title [Construction of free curves by adding lines to a given curve]
{Construction of free curves by adding lines to a given curve}

\author[Alexandru Dimca]{Alexandru Dimca}
\address{Universit\'e C\^ ote d'Azur, CNRS, LJAD, France and Simion Stoilow Institute of Mathematics,
P.O. Box 1-764, RO-014700 Bucharest, Romania}
\email{dimca@unice.fr}

\author[Giovanna Ilardi]{Giovanna Ilardi}
\address{Dipartimento Matematica Ed Applicazioni ``R. Caccioppoli''
Universit\`{a} Degli Studi Di Napoli ``Federico II'' Via Cintia -
Complesso Universitario Di Monte S. Angelo 80126 - Napoli - Italia}
\email{giovanna.ilardi@unina.it}

\author[Piotr Pokora]{Piotr Pokora}
\address{Department of Mathematics,
Pedagogical University of Krakow,
Podchor\c a\.zych 2,
PL-30-084 Krak\'ow, Poland.}
\email{piotr.pokora@up.krakow.pl}

\author[Gabriel Sticlaru]{Gabriel Sticlaru}
\address{Faculty of Mathematics and Informatics,
Ovidius University
Bd. Mamaia 124, 900527 Constanta,
Romania}
\email{gabriel.sticlaru@gmail.com}

\thanks{\vskip0\baselineskip
\vskip-\baselineskip
\noindent A. Dimca was partially supported by the Romanian Ministry of Research and Innovation, CNCS - UEFISCDI, Grant \textbf{PN-III-P4-ID-PCE-2020-0029}, within PNCDI III.\\
P. Pokora was partially supported by the National Science Center (Poland) Sonata Grant Nr \textbf{2018/31/D/ST1/00177}}



\begin{abstract} 
In the present note we construct new families of free  plane curves starting from a curve $C$ and adding high order inflectional tangent lines of $C$, lines joining the singularities of the curve $C$, or lines in the tangent cone of some singularities of $C$. These lines $L$ have in common that the intersection $C \cap L$ consists of a small number of points. We introduce the notion of a supersolvable plane curve and conjecture that such curves are always free, as in the known case of line arrangements. Some evidence for this conjecture is given as well, both in terms of a general result in the case of quasi homogeneous singularities and in terms of specific examples. We construct a new example of maximizing curve in degree 8 and the first and unique known example of maximizing curve in degree 9. In the final section, we use a stronger version of a result due to Schenck, Terao and Yoshinaga  to construct  families of free conic-line arrangements by adding lines to the  conic-line arrangements of maximal Tjurina number recently classified by  V. Beorchia and R. M. Mir\'o-Roig in arXiv:2303.04665.
\end{abstract}
 
\maketitle


\section{Introduction} 

Our goal in this paper is to {\it construct  new free or nearly curves by adding  inflectional tangents or lines passing through the singularities of a given plane projective curve}.  Sometimes, when high order inflectional tangents are missing, lines in the tangent cones of the singularities may replace them successfully. The use of tangent cones is a must when we want to get supersolvable curves with the modular point belonging to a non-linear irreducible component, see Definition \ref{defSS}.

To determine the existence of inflectional tangents, we start by recalling some facts about inflection points.
Let $C:F=0$ be a reduced plane curve in the complex projective plane $\PP^{2}$ which is defined by a homogeneous polynomial $F \in S=\C[x,y,z]$ of degree $d\geq 2$. 
The Hessian of $F$ is given by the following well-known formula
\begin{equation}
\label{eq1}
H=\det \left(
  \begin{array}{ccccccc}
     F_{xx} & F_{xy}& F_{xz}  \\
      F_{xy} & F_{yy}& F_{yz}\\
     F_{xz} & F_{yz}& F_{zz}  \\
   \end{array}
\right).
\end{equation}
Let  $H_C:H=0$ be the Hessian curve associated to $C$.
It is known that the intersection $X_C=C \cap H_C$ consists exactly of the set of inflection points $I_C$ of $C$ union with the set of singular points $Y_C$ of $C$. Recall that
if $p \in C$ is a smooth point of this curve, and $T_pC$ denotes the projective line tangent to $C$ at $p$, then the
{\it inflection order} of $p$ is by definition
\begin{equation}
\label{ip}
\iota_p(C)=(C,T_pC)_p-2,
\end{equation}
where $(C,T_pC)_p$ denotes the intersection multiplicity of the curves
$C$ and $T_pC$ at the common point $p$. Moreover, we say that $p$ is an inflection point of $C$, i.e.,  $p \in I_C$, if and only if $\iota_p(C) >0$.

The  intersection multiplicity $(C,H_C)_p$ of
the curves $C$ and $H_C$ at the point $p$ is a key invariant in understanding the geometry of the plane curve $C$.
When $p\in I_C$ is an inflection point, then the relation between  the inflection order $\iota_p(C)$ of $p$ and the intersection multiplicity $(C,H_C)_p$ is well-known, see for instance \cite[Theorem 9.7 and Corollary 9.10]{K}.
\begin{thm}
\label{thmIO}
For any  reduced plane curve  $C$ of degree $d\geq 2$ and any smooth point $p \in C$, one has 
$$\iota_p(C)=(C,H_C)_p.$$
\end{thm}
It is clear that $\iota_p(C) \leq d-2$, except the case when $p$ sits on a line $L=T_pC$ which is an irreducible component of $C$. In the later case, one has $\iota_p (C)=\infty$, and it is easy to check that the line
$L$ is also an irreducible component of the Hessian curve $H_C$. From now on, {\it while searching for the inflection points, we assume that no irreducible component of $C$ is a line}.

When $p \in Y_C$ is a singular point, then the description of the intersection multiplicity $(C,H_C)_p$ is more subtle. Let $TC_p(C)$ be the reduced projective tangent cone of the curve $C$ at $p$. For any line $L \in TC_p(C)$ we define the corresponding tangential multiplicity
$$m_L(C)=(L,C_L)_p,$$
where $(C_L,p)$ is the union of all branches of the singularity $(C,p)$ whose tangent line at $p$ coincides with $L$. With this notation, one has the following general result, which is a user-friendly reformulation of \cite[Proposition 25]{JP}. 
\begin{thm}
\label{thmA}
Assume that $C$ is a reduced plane curve and that $p \in C$ is any singular point. Then one has
$$(C,H_C)_p=3 \mu_p(C)+m_p(C)-3+\sum_{L \in TC_p(C)}m_L(C),$$
where $\mu_p(C)$ is the Milnor number and $m_p(C)$ is the multiplicity of the singularity $(C,p)$. 
\end{thm}

The case when $(C,p)$ is irreducible is also stated in \cite[Theorem 2.1.9]{Moe}. Indeed, in this case, one has $\mu_p(C)=2\delta_p(C)$ with $\delta_p(C)$ being the $\delta$-invariant of the singularity $(C,p)$. We discuss in detail the statement of Theorem \ref{thmA} and give some examples in the next section.

\medskip
Our new results are the following. Consider  the case when $(C,p)$ is an ordinary $k$-multiple point, that is there are $k$ smooth branches $C_1,\ldots,C_k$ at $p$, with distinct tangent lines $L_1, \ldots, L_k$. If we set $m_j=(C_j,L_j)_p\geq 2$ for $j \in \{1, \ldots, k\}$, then we call $(C,p)$  an ordinary $k$-multiple point of type $(m_1, \ldots, m_k)$. Moreover, we say that $(C,p)$ is an
 {\it ordinary simple $k$-multiple point} if $m_j=2$ for all $j \in \{1, \ldots, k\}$. 
\begin{thm}
\label{thm1}
For any  reduced plane curve  $C$ of degree $d\geq 3$ and any singular point $p \in C$ with multiplicity $m_{p}(C)=k$, one has 
$$(C,H_C)_p \geq 3k(k-1),$$
and equality holds if and only if $(C,p)$ is an ordinary simple $k$-multiple point.
\end{thm}

\begin{cor}
\label{cor1}
For any  reduced plane curve  $C$ of degree $d\geq 3$ and having $s$ singular points, one has
$$i(C) \leq 3d(d-2)-6s.$$
More precisely, if $n_k$ denotes the number of singular points of $C$ of multiplicity $k$, then
$$i(C) \leq 3d(d-2)-\sum_k 3k(k-1)n_k.$$
The equality occurs if and only if all the singularities of $C$ are ordinary simple $k$-multiple points for various $k$. 
\end{cor}

This result is a significant improvement of the inequality
$$\sum_{p \in I_C} \iota_p \leq 3d(d-2)-s$$
which is given in \cite[Corollary 9.9]{K}. Our next result is a general construction of free curves by adding inflectional tangents and lines passing through the singularities of the initial curve $C$, which is a curve of Thom-Sebastiani type. Let $m \geq 2$ be a positive integer and let $\ell_j(x,y)$  with $j \in \{1, \ldots, m\}$ be $m$ distinct linear forms in $x$ and $y$. Consider the curve
$$C: F=\prod_{j=1}^m\ell_j(x,y)^{k_j} -z^{d}=0,$$
in $\PP^2$, where $k_j \geq 1$ are positive integers such that $\sum k_j=d$, 
and the family of lines $L_j:\ell_j(x,y)=0$. The line $L_j$  is the inflectional tangent at the point $p_j=(x_j:y_j:0)$ where $(x_j:y_j) \in \PP^1$ is the zero set of $\ell_j$ whenever $k_j=1$, and it is the reduced tangent cone at the singularity $p_j$ when $k_j>1$. 

\begin{thm}
\label{thm2}
With the notation as above, the curve
$$C'=C \cup \bigcup_{j=1}^mL_j : \quad F' =F \cdot \prod_{j=1}^m\ell_j=0$$
of degree $d+m$ is free with the exponents $(m-1,d)$. Moreover, if $L:z=0$ is the line passing through all the points $p_j$ of $C$,
then  the curve
$$C''=C' \cup L: \quad F'' =zF \cdot\prod_{j=1}^m\ell_j=0$$
of degree $d+m+1$ is free with the exponents $(m-1,d+1)$.
\end{thm}
Starting with a smooth Fermat type curve and adding all of its inflectional tangents and the $3$ coordinate axes, one gets again a free curve, as the following result shows. 
Consider the Fermat curve $C:x^d+y^d+z^d=0$. Let $\epsilon$ be any root of the equation $t^d+1=0$ in $\C$. Then the line
$L_{\epsilon}: y=\epsilon x$ intersects $C$ only at the point
$p_{\epsilon}=(1:\epsilon: 0)$. Hence $p_{\epsilon}$ is an inflection point of the maximal order and hence $i(L_{\epsilon})=d-2$, i.e., one has the equality in \eqref{eq3.1}. 
In this way, we get $d$ inflection points, which are just the intersection of $C$ with the line $z=0$. 
Cyclically permuting
$x,y,z$, we get all the $3d$ inflection points of this type, and these are all the inflection points of $C$ by \eqref{eq2.1}. The following result is a special case of Theorem \ref{thm2}, when $k_j=1$ for all $j$.

\begin{cor}
\label{corT2}
The union $C': F'=(x^d+y^d)F=0$ of the smooth Fermat curve $C: F=x^d+y^d+z^d=0$
of degree $d$ with the $d$ inflectional tangents $L_{\epsilon}$ meeting at one point, is a free curve of degree $2d$ and the exponents are $(d-1,d)$.
When $d=3$, the curve $C'$ is maximizing of degree $6$.
The union $C'': F'' =zF'=0$ of the curve $C'$ with the line  $L:z=0$  passing through all the flex points $p_{\epsilon}$ of $C$,
is a free curve of degree $2d+1$ and the exponents are $(d-1,d)$.
When $d=3$, then the curve $C'$ is maximizing of degree $7$.
\end{cor}
Note that if we continue to add just inflectional tangents, the resulting curve is no longer a free curve. For instance, the curve
$$C'': F''=(x^3+y^3+z^3)(x^3+y^3)(y^3+z^3)=0$$
is nearly free with the exponents $(4,5)$, and the curve
$$C''': F'''=(x^3+y^3+z^3)(x^3+y^3)(y^3+z^3)(x^3+z^3)=0$$
is not even nearly free.  However, we have the following general result.
\begin{thm}
\label{thm3}
The smooth Fermat curve
$C: F=x^d+y^d+z^d=0$
 has exactly $3d$ inflectional tangents and their union forms the following line arrangement
$$\A: (x^d+y^d)(y^d+z^d)(x^d+z^d)=0.$$
The union of the curve $C$, the lines in $\A$, and the $3$ coordinate axes produce a new curve
$$C': F'=xyz (x^d+y^d)(y^d+z^d)(x^d+z^d)(x^d+y^d+z^d) =0$$
of degree $4d+3$, which is free with the exponents $(2d+1,2d+1)$. 
Moreover, the curve $C'' \subset C'$ given by
$$C'':F''=xy(y^d+z^d)(x^d+z^d)(x^d+y^d+z^d) =0$$
has degree $3d+2$ and it is free with the exponents $(d+1,2d)$. \\
When $d=2$, then the curve $C''$ is maximizing of degree $8$.
\end{thm}
It is easy to prove that  the curve
$$C: F=x^my^m+y^mz^m+x^mz^m=0$$
has no inflection points using Theorem \ref{thmA}, see Example \ref{ex3} below. To get a free curve from $C$, we may add two of the three tangent cones, or just one tangent cone and two lines joining singular points. Indeed, one has the following result.

\begin{thm}
\label{thm4}
The curve
$$C': F'=(x^my^m+y^mz^m+x^mz^m)(x^m+y^m)(y^m+z^m)=0$$
has degree $4m$ and it is free with exponents $(2m-1,2m)$ for any $m \geq 2$. The curve 
$$C'': F''=yz(x^my^m+y^mz^m+x^mz^m)(y^m+z^m)=0$$
has degree $3m+2$ and it is free with exponents $(m+1,2m)$ for any $m \geq 2$.
\end{thm}
\begin{definition}
\label{defSS}
Given a reduced plane curve $C$, we say that $p \in C$ is a {\it modular point} for $C$ if the central projection 
$$\pi_p : \PP^2 \setminus \{p\} \to \PP^1$$
induces a locally trivial fibration of the complement $M(C)=\PP^2 \setminus C$. We say that a 
reduced plane curve $C$ is {\it supersolvable} if it has at least one modular point.
\end{definition}
The map induced by $\pi_p$ is a locally trivial fibration if and only if for any line $L_p$ passing through $p$ and not an irreducible component of $C$, one has
$$(C,L_p)_p=\mult_p(C) \text{ and } (C,L_p)_q=1 \text{ for any  } q \in C \cap L_p, \ q \ne p.$$
This fibration has as base and as fiber a projective line
$\PP^1$ with a number of points deleted, and hence the complement $M(C)$ is a $K(\pi,1)$ space. When $C$ is a line arrangement, this definition of a modular point coincides with the usual one,  and a line arrangement is supersolvable by definition if it has a modular point. In particular, the existence of a modular point for a line arrangement $C$ implies that $C$ is free,  see for all these well known facts \cite{HA,OT}.  We venture to make the following.
\begin{conj}
\label{conj1}
A supersolvable plane curve $C$ is free.
\end{conj}
One setting where this conjecture holds is the following.

\begin{thm}
\label{thm5}
Let $C$ be a reduced plane curve, let $p \in M(C)$ be a point and let $\A$ be the set of lines $L$ passing through $p$ such that there is a point $q \in L \cap C$ with $(C,L)_q >1$. Assume that all the singularities of the curve $C'$ obtained by adding all the lines in $\A$ to $C$ are quasi homogeneous. Then $C'$ is supersolvable and free.
In particular, this holds when all the singularities $s$ of $C$ have multiplicity 2, and $p$ is not on any tangent cone $TC_s(C)$ for $(C,s)$ a singularity with $\mu(C,s) \geq 3$.
\end{thm}
When $C$ is a plane curve having only nodes $A_1$ and cusps $A_2$ as singularities, and in addition $p$ is a generic point, this result is known, see \cite[Theorem 1.12]{Mich}. Moreover, it is easy to see that the free curves $C'$ and $C''$ constructed above  in Theorem \ref{thm2}
are special cases of the construction in Theorem \ref{thm5}.
On the other hand, the  free curve $C''$ constructed above  in Theorem \ref{thm4} is of a different nature, since in this case $p \in C$ and the resulting curve $C''$  has not only quasi homogeneous singularities. However, this curve gives new examples where Conjecture \ref{conj1} holds, in view of the following result.

\begin{prop}
\label{prop5}
The free curve  $C''$  constructed in Theorem \ref{thm4} is supersolvable. In particular,  the associated complement $M(C'')$ is a $K(\pi,1)$ space.
\end{prop}
We explain in Example \ref{exKP} that the other free curves constructed above in Theorems \ref{thm3} and \ref{thm4} are not supersolvable.

The organization of the paper goes as follows. In Section $2$, we discuss the proof of Theorem \ref{thmA} explaining all the necessary details. In Section $3$ we recall basic facts on the free, nearly free and maximizing curves.

In Section 4 we deliver our proofs of Theorems \ref{thm1}, \ref{thm2},  \ref{thm3}, \ref{thm4} and  \ref{thm5}, and of Proposition \ref{prop5}. Then, in Section $5$, we describe all smooth plane quartic curves admitting the maximal possible number of flexes of maximal order, i.e., flexes of order $2$. There are two such curves, and only one of them, the Fermat curve, yields  free curves as in Corollary \ref{corT2} and in Theorem \ref{thm3}. In Section $6$ we discuss some singular plane curves and the free curves obtained from them, which are sometimes supersolvable as well, given rise to supersolvable  free curves not covered by our general Theorem \ref{thm5}. In Example \ref{ex4} we construct a {\it new example of maximizing curve in degree 8 and the first and unique known example of maximizing curve in degree} 9.
In the final section, we use a stronger version of a result due to Schenck, Terao and Yoshinaga in \cite{STY} to construct {\it families of free conic-line arrangements} starting with the free conic-line arrangements
$C:F=0$ with exponents $(1,d-2)$, where $d= \deg F$, which have been recently classified by  V. Beorchia and R. M. Mir\'o-Roig in \cite{Be}. These families contain free conic-line arrangements with arbitrary exponents and also provide countable examples where the Conjecture \ref{conj1} holds even in the presence of non quasi homogeneous singularities, see Remarks \ref{rkDIS1} and \ref{rkDIS2}.

\bigskip

We would like to thank the referee for his/her careful reading of the manuscript and the useful suggestions.

\section{Discussion on Theorem \ref{thmA} and some examples} 
The paper \cite{JP} uses rather heavy notations, and perhaps due to this fact has a smaller impact than it deserves. Let us introduce some notation. For a reduced plane curve $C:F=0$ and any point $q=(\al:\be: \gamma) \in \PP^2$, we define the polar $\Delta_q(C)$ of $C$ with respect to $q$ by the equation
\begin{equation}
\label{eq4}
\Delta_q(C): \al F_x+ \be F_y + \gamma F_z =0.
\end{equation}
For a property $\Pc$, the authors of \cite{JP} use the notation $\bf 1 _{\Pc}$ to denote $1$ if $\Pc$ is true and 0 otherwise, see the discussion before Theorem 2 in \cite{JP}. The first equality in \cite[Proposition 25]{JP} gives the expression of the intersection multiplicity
$(C,\Delta_q(C))_p$ for any singular point $p \in C$ and any point $q \in \PP^2$. Using our discussion above, we see that for $q\ne p$, $q$ not on any line $L$ in the tangent cone $TC_p(C)$, this multiplicity $(C,\Delta_q(C))_p$ is given by a double sum $\SSS$, i.e., the
second sum involving the characteristic functions $\bf 1 _{\Pc}$ vanishes. With this observation, the second equality in \cite[Proposition 25]{JP} can be stated as
\begin{equation}
\label{eq5}
(C,H_C)_p=3(C,\Delta_q(C))_p+I_p,
\end{equation}
where $I_p=\sum_{i \in \I}(i^{(i)}_p-2)$. Here $\I$ is a set of indices parametrizing the pro-branches of $(C,p)$ and $i^{(i)}_p$ is the tangential intersection number of the pro-branch $B _i$, see \cite[Definition 22]{JP}. Recall that any branch $B$ of a plane curve singularity $(C,0)$ at the origin of $\C^2$, such that $x=0$ is not a tangent line, can be defined by a Weierstrass polynomial
$$\Gamma_B(x,y)=\prod_{j=1,m_B}(y-\phi_{B,j}(x))$$
where $m_B$ is the multiplicity of the branch $B$ and there is an analytic function $\phi_B(x) \in \C\{x\}$ with $\ord \phi_B(x) \geq m_B$ such that
\begin{equation}
\label{eq5R}
\phi_{B,j}(x)=\phi_B(\exp(2\pi i j/m_B)x^{\frac{1}{m_B}}).
\end{equation}
With this notation, the branch $B$ has $m_B$ associated pro-branches
$$B_j: y-\phi_{B,j}(x)=0$$
and the corresponding tangential intersection number is given by
$$i^{B_j}_0=\ord (\phi_{B,j}(x) -\phi_{B,j}'(0)x) \in \Q.$$
It follows from equation \eqref{eq5R} that $i^{B_j}_0=i^{B_k}_0$ for any $0 \leq j \leq k \leq m_B$
and hence 
$$m_Bi^{B_j}_0=(L,C_L)_0,$$
for any $j$, with $L$ the tangent line to $B$ and $C_L=B$. This discussion implies that one has
\begin{equation}
\label{eq6}
I_p=\sum_{i \in \I}(i^{(i)}_p-2)=\sum_{L \in TC_p(C)}m_L(C)-2m_p(C),
\end{equation}
since clearly $|\I|=m_p(C)$, each branch having exactly a number of pro-branches given by the multiplicity of that branch. 
Next, we return to the intersection multiplicity  $(C,\Delta_q(C))_p$. We can assume that $p=(0:0:1)$ and set $f(x,y)=F(x,y,1)$, then one has
$$f_x(x,y)=F_x(x,y,1), \ f_y(x,y)=F_y(x,y,1) \text{, and } xf_x+yf_y+F_z(x,y,1)=f.$$
If we define the generic local polar variety of the singularity
$$(C,0):f(x,y)=0$$
by the equation
$$\Delta_0(C): \al 'f_x+ \be' f_y  =0,$$
with $(\al': \be') \in \PP^1$  being a generic point, it is easy to see that
\begin{equation}
\label{eq6.1}
(C,\Delta_0(C))_0=(C,\Delta_q(C))_p.
\end{equation}
In fact, the line $L':z=0$ is clearly not in the tangent cone $TC_p(C)$ since $p \notin L'$, and hence we may take $\gamma=0$ 
and $(\al: \be) \in \PP^1$  generic  in the formula  \eqref{eq4}.
In order to compute this local intersection number
$$\kappa_0(C)=(C,\Delta_0(C))_0,$$
which is also called the $\kappa$-invariant of the singularity $(C,0)$, we can use for instance \cite[Proposition 3.38]{GLS} and get
\begin{equation}
\label{eq7}
\kappa_0(C)=\mu(C,0)+m_0(C,0)-1.
\end{equation}
If we use this formula for the singularity $(C,p)$, we get from
\eqref{eq5} and \eqref{eq6} the following equality
$$(C,H_C)_p=3(\mu_p(C)+m_p(C)-1)+\sum_{L \in TC_p(C)}m_L(C)-2m_p(C)=$$
$$=3 \mu_p(C)+m_p(C)-3+\sum_{L \in TC_p(C)}m_L(C).$$
This proves our reformulation of \cite[Proposition 25]{JP} in Theorem \ref{thmA}.

\medskip

In order to construct free plane curves starting with a curve $C$ of degree $d$, by adding lines, in particular inflectional tangents,
we have to look for lines
$L$ such that the sum
\begin{equation}
\label{eq2}
i(L)=\sum _{p \in L \cap I_C, \  T_pC=L}\iota_p(C)
\end{equation}
is as large as possible with respect to the degree $d$. Given a curve $C$, first we use Theorem \ref{thmA} to count the total number of inflection points of $C$, namely
\begin{equation}
\label{eq2.1}
i(C)=\sum _{p \in I_C}  \iota_p(C)=3d(d-2)-\sum _{p \in Y_C} (C,H_C)_p.
\end{equation}
One clearly has for any line $L$
\begin{equation}
\label{eq3}
i(L) \leq i(C),
\end{equation}
and the equality holds if and only if for any point $p \in I_C$ one has $T_pC=L$. Moreover,
\begin{equation}
\label{eq3.1}
i(L) \leq \sum _{p \in L \cap C}((C,L)_p-2)=d-2|L\cap C|,
\end{equation}
and the equality holds if and only if for any point $p \in L \cap C$ one has  $T_pC=L$.

\begin{ex}
\label{exA1}
Let us consider the case when $(C,p)$ is a node $A_1$, that is there are two smooth branches $(C_1,p)$ and $(C_2,p)$ meeting transversally at $p$. Let
$T_1=T_pC_1$ and $T_2=T_pC_2$ be the associated tangent lines and define {\it the type of the node } $(C,p)$ to be the pair of integers
$$(m_1,m_2)=((C_1,T_1)_p, (C_2,T_2)_p).$$
It is clear that $m_j \geq 2$ for $j=1,2$. When $m_1=m_2=2$, then $(C,p)$ is said to be a {\it simple node}, and one knows that $(C,H_C)_p=6$, see \cite[pp. 68--69]{F}. In the general situation, Theorem \ref{thmA} gives the equality
$$(C,H_C)_p=3+2-3+m_1+m_2=m_1+m_2+2.$$
More generally, consider  the case when $(C,p)$ is an ordinary $m$-multiple point, that is there are $m$ smooth branches $C_1,\ldots,C_m$ with distinct tangent lines $L_1, \ldots, L_m$. If we set $m_j=(C_j,L_j)_p$ for $j=1, \ldots,m$, then 
  Theorem \ref{thmA} gives the equality
$$(C,H_C)_p=3(m-1)^2+m-3+\sum_{j=1,m}m_j=m(3m-5)+\sum_{j=1,m}m_j.$$
\end{ex}

\begin{ex}
\label{exA2} Let us consider the case when $(C,p)$ is a
 singularity $A_{2m-1}$ with  $m\geq 2$. Then there are two tangent smooth branches with a common tangent line $L$.  If we set $m_L=(C,L)_p$, then Theorem \ref{thmA} gives the equality
$$(C,H_C)_p=3(2m-1)+2-3+m_L=2(3m-2)+m_L.$$

\end{ex}

\begin{ex}
\label{exA3} Let us consider the case when $(C,p)$ is a
 singularity $A_{2m}$ with  $m\geq 1$. Then there is a unique  branch, with a  tangent line $L$.  If we set $m_L=(C,L)_p$, then Theorem \ref{thmA} gives the equality
$$(C,H_C)_p=3(2m)+2-3+m_L=6m-1+m_L.$$
When $(C,p)$ is a
 cusp $A_{2}$, only the value $m_L=3$ is possible, and hence
 $(C,H_C)_p=8$ in this case. 
\end{ex}

\section{Free, nearly free and maximizing curves} 

In this section we recall some basic facts on free, nearly free and, maximizing curves in $\PP^2$  following \cite{MaxCurv,FNF}.

Let 
$${\rm Der}(S) = \{ \partial := a\cdot \partial_{x} + b\cdot \partial_{y} + c\cdot \partial_{z}, \,\, a,b,c \in S\}$$
 be the free $S$-module of $\mathbb{C}$-linear derivations of the polynomial ring $S$. For a reduced curve $C \, : F=0$, we introduce
$${\rm D}(F) = \{ \partial \in {\rm Der}(S) \, : \, \partial\, F \in \langle F \rangle \},$$
the graded $S$-module of derivations preserving the ideal $\langle F \rangle$. We have the following decomposition
$${\rm D}(F) = {\rm D}_{0}(F) \oplus S\cdot \delta_{E},$$
where $\delta_{E} = x\partial_{x} + y\partial_{y} + z\partial_{z}$ is the Euler derivation and 
$${\rm D}_{0}(F) = \{ \partial \in {\rm Der}(S) \, : \, \partial \, F = 0\}$$
is the set of all $\mathbb{C}$-linear derivations of $S$ killing the polynomial $F$.
\begin{definition}
We say that a reduced curve $C \, : F=0$  is \textbf{free} if ${\rm D}(F)$, or equivalently ${\rm D}_{0}(F)$, is a free graded $S$-module. The exponents $(d_1,d_2)$ of a  free curve $C$  are the degrees of a basis for the free graded $S$-module ${\rm D}_{0}(F)$ which rank 2.
\end{definition}
\begin{rk}
\label{rkF}
The exponents $(d_1,d_2)$ of a free curve $C:F=0$ of degree $d$ are known to satisfy $d_1+d_2=d-1$. Conversely, if there are two elements
$r_1,r_2 \in {\rm D}_{0}(F)$, which are $S$-linearly independent and satisfy
$$d_1+d_2=d-1$$
then the curve $C$ is free with exponents $(d_1,d_2)$, see \cite{ST,T}.
\end{rk}

\begin{definition}
The minimal degree of derivations killing $F$, or of Jacobian syzygies involving the partial derivatives of $F$, is defined as
$${\rm mdr}(F) = {\rm min}\{r \in \mathbb{N} \, : \, {\rm D}_{0}(F)_{r} \neq 0 \}.$$
\end{definition}
To check  whether a given plane curve is free, one may use the following result by du Plessis and Wall \cite{duP}.

\begin{thm}
\label{thmCTC}
Let $C: \,F=0$ be a reduced plane curve of degree $d$, let $r = {\rm mdr}(F)$ and let $\tau(C)$ be the total Tjurina number of $C$.
Then the following two cases hold.
\begin{enumerate}
\item[a)] If $r < d/2$, then $\tau(C) \leq \tau(d,r)_{max}= (d-1)^2-r(d-r-1)$ and the equality holds if and only if the curve $C$ is free.
\item[b)] If $d/2 \leq r \leq d-1$, then
$\tau(C) \leq \tau(d,r)_{max}'$,
where, in this case, we set
$$\tau(d,r)_{max}'= \tau(d,r)_{max} - \binom{2r-d+2}{2}.$$
\end{enumerate}
\end{thm}
\begin{definition}
A reduced curve $C: F=0$ of degree $d$ is \textbf{nearly free} if either $\mdr(F) < d/2$ and  $\tau(C) =\tau(d,r)_{max}-1$, or $ \mdr(F) = d/2$ and  $\tau(C) =\tau(d,r)_{max}$. In addition, the exponents of a nearly free curve $C:F=0$ of degree $d$ are given by the pair $(\mdr F, d-\mdr F)$.
\end{definition}

\begin{definition}
\label{defMC}
A curve $C:F=0$ of degree $d$ having only ${\rm ADE}$-singularities is \textbf{maximizing} if either 
$d=2m$ and  $\tau(C) =3m(m-1)+1$, or  $d=2m+1$ and  $\tau(C) =3m^2+1$.
\end{definition}
The relation between maximizing curves and free curves is the following, see \cite{MaxCurv}.

\begin{thm}
\label{thmMC}
A curve $C:F=0$ of degree $d$ having only ${\rm ADE}$-singularities is {maximizing} if and only if either 
$d=2m$ and $C$ is a free curve with the exponents $(m-1,m)$, or  $d=2m+1$ and $C$ is a free curve with the exponents $(m-1,m+1)$.
\end{thm}
In the sequel we need the following version of \cite[Theorem 1.10]{Mich}.
\begin{thm}
\label{thmMich}
Let $C: \,F=0$ be a reduced plane curve of degree $d$ and let $p$ be any point of $C$. Let $\A$ be the union of the irreducible components of $C$ which are lines passing through $p$, and let $C':G=0$ be the union of the other irreducible components of $C$. We assume that $p \in C'$. Let $m=|\A|$ and $e=\deg G$. Then $r = {\rm mdr}(F)$ can be in one of the following cases.
\begin{enumerate}
\item[a)] $r=e$;
\item[b)] $r=m-1$ and $C$ is free with exponents $(m-1,e)$;

\item[c)] $m \leq r \leq e-1$.

\end{enumerate}
\end{thm}
The only difference of this result with respect to  \cite[Theorem 1.10]{Mich} is that here $p \in C'$.
\proof
The key part of the proof of \cite[Theorem 1.10]{Mich} is contained in
\cite[Lemma 4.3]{Mich}, where a Jacobian syzygy
$$\rho: aF_x+bF_y+cF_z=0$$
is constructed using a differential 2-form $\omega$, and it is shown that this syzygy is primitive, that is there is no common factor for $a,b,c\in S$.
It is in this latter part that the condition $p \notin C'$ was used.
If we assume that $p=(0:0:1)$, it is easy to see that the Jacobian relation 
$\rho$ constructed there is still a Jacobian syzygy in our situation. Moreover, it is a primitive syzygy if and only if $G$ and $G_z$ have no common factor.
Let $M\in S$ be an irreducible polynomial which is a common factor for 
$G$ and $G_z$. Note that $M$ cannot involve only $x$ and $y$, since this would correspond to a line in $C'$ passing through $p$, which is impossible by the definition of $C'$. It follows that $M_z \ne 0$.
If $G=MN$, then $G_z=M_zN+MN_z$, which implies that either
$M$ divides $N$ or $M$ divides $M_z$. But $M$ cannot divide $N$, since $C'$ is reduced. And $M$ cannot divide $M_z$, since the degree of $M_z$ as a polynomial in $z$ is strictly smaller than the corresponding degree of $M$. This contradiction proves our claim.
\endproof

\section{The proofs of our main results} 
\subsection{Proof of Theorem \ref{thm1}}
We can assume that $p=(0:0:1)$ and set $f(x,y)=F(x,y,1)$, $g(x,y)=j^kf(x,y)$ the initial form of $f$, that is the sum of the lowest degree terms in the Taylor expansion of $f$ at $0$. The notation $j^kf(x,y)$, the $k$-th jet of $f$ at $0$, is an alternative way of notation for this binary form of degree $k$. Let $n$ be the number of distinct factors of $g$. If $n=k$, then using Example \ref{exA1} we have
$$(C,H_{C})_p=k(3k-5)+\sum_{j=1,k}m_j \geq 3k(k-1),$$
since $m_j \geq 2$. Moreover, it is obvious that the equality holds if and only if $(C,p)$ is an ordinary simple $k$-multiple point.

Assume now that $n<k$ and let
$$g(x,y)=\ell_1^{a_1}\cdot ... \cdot \ell_n^{a_n}$$
be the decomposition of $g$ as a product of linear factors. Recall that for two isolated plane curve singularities $(X,0)$ and $(Y,0)$ with no common component one has 
\begin{equation}
\label{eqMilnor}
\mu(X\cup Y,0)=\mu(X,0)+\mu(Y,0)+2(X,Y)_0-1,
\end{equation}
see \cite[Theorem 6.5.1]{CTC}.
Let $C_j:f_j=0$ be the union of the branches of $(C,p)$ which are tangent to the line $L_j: \ell_j=0$ with $j \in \{1, \ldots, n\}$.
Then $(C,p)=(C_1,p) \cup ... \cup (C_n,p)$ and we estimate $\mu_p(C)=\mu(C,p)$ using the above formula. To start with, note that
since $j^{a_j}f_j=\ell_j^{a_j}$, it follows that
$$\mu(C_j,p) \geq a_j(a_j-1).$$
We prove by induction on $m$ that 
$$\mu((C_1,p) \cup ... \cup (C_m,p))\geq (a_1+ \ldots + a_m)^2-(a_1+ \ldots + a_m)-m+1$$
for any $1 \leq m \leq n$. This inequality holds for $m=1$ as we have already seen above. Assume that the inequality holds for some $m<n$. 
Then it follows that
$$\mu((C_1,p) \cup ... \cup (C_m,p)\cup (C_{m+1},p))\geq ((a_1+ \ldots + a_m)^2-(a_1+ \ldots + a_m)-m+1)+$$
$$+a_{m+1}(a_{m+1}-1)+2a_{m+1}(a_1+ \ldots + a_m)-1=$$
$$=(a_1+ \ldots + a_{m+1})^2-(a_1+ \ldots + a_{m+1})-m,$$
which completes our proof by induction.
Since $a_1+ \ldots + a_n=k$, this yields the inequality
$$\mu(C,p) \geq k^2-k-n+1.$$
On the other hand, one has
$$m_{L_j}=(C_j,L_j)_p \geq a_j+1$$
and hence
$$\sum_jm_{L_j} \geq k+n.$$
Using Theorem \ref{thmA}, we get
$$(C,H_{C})_p \geq 3(k^2-k-n+1)+k-3+k+n=3k(k-1)+2(k-n)>3k(k-1),$$
since we have assumed $n<k$.
This completes the proof of Theorem \ref{thm1}. 

Corollary \ref{cor1} is an obvious consequence of Theorem \ref{thm1}.

\subsection{Proof of Theorem \ref{thm2}}
The curve $C$ was studied in \cite[Example 4.5]{3-syz} and it was shown that the minimal degree of a Jacobian relation for $C$ is given by
$$\mdr(F)=m-1.$$
Since the minimal degree of a Jacobian relation can only increase when one adds lines to a given curve, see \cite[Proposition 3.1]{mdr}, it follows that 
$$\mdr(F') \geq \mdr(F)=m-1.$$
The curve $C$ has $m$ inflection points and singularities on the line $z=0$, locally given by equations $u^{k_j}-v^d=0$, located at the points $p_j=(a_j:b_j:0)$, with $\ell_j(a_j,b_j)=0$ for $j \in \{1,\ldots, m\}$.
When we add the line $L_j$, we get at the point $p_j$ a weighted homogeneous singularity of degree $d_j=1$ with respect to the weights
$w_1=wt(u_j)=\frac{1}{k_j+1}$ and $w_2=wt(z)=\frac{k_j}{d(k_j+1)}$, where $u_j=\ell_j$ is a local coordinate at $p_j$ on the line $L_j$.
It follows the following equality involving Tjurina and Milnor numbers:
$$\tau(C',p_j)=\mu(C',p_j)=\frac{(1-w_1)(1-w_2)}{w_1w_2}=(d-1)k_j+d.$$
Hence the total Tjurina number of $C'$ is
$$\tau(C')=d(d-1)+md+(m-1)^2,$$
since clearly 
$$\tau(C',p)=\mu(C',p)=(m-1)^2.$$
Now a curve of degree $d'=d+m$ with $r'=\mdr(F') $ satisfies the inequality
$$\tau(C') \leq \tau(d',r')_{max}$$
where the function 
$$\tau(d',r')_{max}=(d'-1)^2-r'(d'-1-r')$$ 
is a decreasing function of $r'$ for $2r'<d'$, which follows from Theorem \ref{thmCTC}, and the equality $\tau(C') =\tau(d',r')_{max}$ implies that $2r'<d'$ and $C'$ is free with the exponents $(r',d'-r'-1)$. In our case, we get
$$\tau(d',m-1)_{max}=(d+m-1)^2-d(m-1)=\tau(C'),$$
and this proves our claim. The proof of the second claim goes analogously.
\begin{rk}
\label{rkTS}
One can check that the curve $C'$, resp. $C''$, can be regarded as a special case of the curve $C$ constructed in Theorem \ref{thm5}, starting from the curve
$C_0:F=0$ and $p=(0:0:1)$ for $C'$, resp. $C_0:zF=0$ and $p=(0:0:1)$ for $C''$. This is an alternative way to proving Theorem \ref{thm2}.
\end{rk}
\subsection{Proof of Theorem \ref{thm3}}

We start with the following.
\begin{lem}
\label{lem1}
Consider the line arrangement
$$\B: g=xyz (x^d+y^d)(y^d+z^d)(x^d+z^d)=0.$$
Then $\mdr(g)=2d+1$.
\end{lem}
\proof
We consider first the subarrangement of $\B$ given by
$$\B_0: g_0=xyz (x^d+y^d)(y^d+z^d)=0.$$
Note that in this arrangement $\B_0$ there are two points
of multiplicity $d+2$, connected by the line $y=0$. All the lines pass through one of these two points, and the other intersection points are all double points. It follows that $\B_0$ is a line arrangement of type $\hat L(d+2,d+2)$, as in \cite[Definition 4.9]{DIM}, and 
$$\mdr(g_0)=d+1,$$
see \cite[Example 4.11]{DIM}. To get the arrangement $\B$ from $\B_0$,
we have to add  the $d$ lines $L_1,...,L_d$ given by $x^d+z^d=0$. At the stage $k$, where $1 \leq k \leq d$, we have to add the line $L_k$ to the arrangement
$$\B_k=\B_0 \cup L_1 \cup ... \cup L_{k-1}.$$
Note that the intersection of $L_k$ and $\B_k$ consists of exactly $2d+2$ points. If $\B_k$ is given by the reduced equation $f_k=0$, for 
$1 \leq k \leq d$, it follows from \cite[Corollary 6.4]{mdr} that one has
$$\mdr(f_k)=\mdr(f_{k-1})+1$$
for all $1 \leq k \leq d$. Hence 
$$\mdr(g)=\mdr(f_d)=d+1+d=2d+1.$$
\endproof

Using \cite[Theorem 5.1 (b)]{mdr} we see that ${\rm mdr}(F')\geq 2d+1$.
We compute now the total Tjurina number of the curve $C'$. This curve has $3d^2$ nodes $A_1$ and 3 points with local equation $uv(u^d+v^d)$ coming from the double points of the line arrangement $\B$ not situated on $xyz=0$ and the 3 points of multiplicity $d+2$.
The line $x=0$ in $\B$ contains $d$ double points of this line arrangement, which are precisely the inflection points of order $d-2$ of the Fermat curve $C$ situated on this line. The corresponding inflectional tangents are the lines given by $y^d+z^d=0$.
Therefore, when we add $C$, each of these $d$ points becomes a singularity of type $D_{2d+2}$. Similar remarks apply to the lines $y=0$ and $z=0$. It follows that
$$\tau(C')=3d^2+3(d+1)^2+3d(2d+2)=12d^2+12d+3.$$
On the other hand, we have 
$$\tau(4d+3,2d+1)_{max}=(4d+2)^2-(2d+1)^2=12d^2+12d+3.$$
The equality 
$$\tau(C')=\tau(4d+3,2d+1)_{max}$$
implies as above that $r'=\mdr(F')=2d+1$ and that $C'$ is a free curve with the exponents $(2d+1,2d+1)$. The claims for the curve $C''$ are proved in a similar way. The line arrangement
$$\B': g'=xy(y^d+z^d)(x^d+z^d)=0$$
satisfies $\mdr(g')=d+1$, see \cite[Proposition 4.10]{DIM}.
The lines $x=0$ and $y=0$ contain each $d$ points of type $D_{2d+2}$
as above. Besides these points, the line arrangement $\B'$ has two points of multiplicity $d+1$ and $d^2+1$ double points. It follows that
$$\tau(C'')=2d(2d+2)+2d^2+(d^2+1)=7d^2+4d+1=\tau(3d+2,d+1)_{max}.$$

\begin{rk}
\label{rkT3} The line arrangement $\B$ considered in Lemma \ref{lem1} is clearly a subarrangement of the line arrangement
$$\CC: xyz(x^{2d}-y^{2d})(y^{2d}-z^{2d})(x^{2d}-z^{2d})=0,$$
and hence $\B$ is a {\it triangular arrangement} as defined in \cite{MaVa}. One can obtain an alternative proof of Lemma \ref{lem1} using results from Section 4 in \cite{MaVa}. 
 It is interesting to note that the line arrangement formed by the corresponding $3d$ inflectional tangent lines is  the arrangement
$$\B'': (x^d+y^d)(y^d+z^d)(x^d+z^d)=0,$$
which is far from being free. This can be seen using \cite[Theorem 5.1]{MaVa}, since $\B''$ is itself a triangular arrangement.
On the other hand, the line dual to the inflection point $p_{\epsilon}=(1:\epsilon: 0)$ is $L_{\epsilon}':x+\epsilon y=0$, and the union of all these $d$ dual lines obtained when $\epsilon$ varies, is given by
$x^d-y^d=0$ when $d$ is odd. Therefore, for $d$ odd, the line arrangement formed by the corresponding $3d$ dual lines is precisely the free monomial (or Fermat) line arrangement
$$(x^d-y^d)(y^d-z^d)(x^d-z^d)=0.$$
The case $d=3$ is, of course, well-known.

\end{rk}

\subsection{Proof of Theorem \ref{thm4}}
First we consider the curve $C'$.
The reader can check the following Jacobian relations $r_1,r_2 \in {\rm D}_0(F')$ 
$$r_1:  z^{m-1}(x^m+y^m)F'_x-x^{m-1}(y^m+z^m)F'_z=0$$
and
$$r_2: xy^{m-1}(2x^m+3y^m)F'_x-(2F+y^{2m})F'_y+zy^{m-1}(2z^m+3y^m)F'_z=0,$$
where $F=x^my^m+y^mz^m+x^mz^m$. Since 
$$\deg r_1 + \deg r_2= (2m-1)+2m=\deg F' -1,$$
our claim is proved by Remark \ref{rkF}.

We consider now the curve $C''$. First we show that $\mdr F''=m+1$.
To do this, we first determine a minimal degree Jacobian syzygy $r_1$ for $F$.
One has 
$$r_1: a_1F_x+b_1F_y+c_1F_z=0,$$
where $a_1= x(y^m-z^m)$, $b_1=-y(y^m+z^m)$ and $c_1=z(y^m+z^m)$. Now we apply \cite[Theorem 3.3]{mdr} and see that if we add a line $L_0$ 
to $C$ given by an equation $\ell \, : \,sy+ty=0$ such that $\ell$ divides
$$sb_1+tc_1=(-sy+tz)(y^m+z^m),$$
the resulting curve $C_0=C \cup L_0: F_0=\ell F=0$ has again
$$\mdr F_0= \mdr F=m+1.$$
Moreover, the coefficients of a minimal degree syzygy for $F_0$ can be obtained from the discussion just before  \cite[Theorem 3.3]{mdr}. It follows that one can add one by one all the lines in the arrangement
$$yz(y^m+z^m)=0$$
and get at the end $\mdr F''=\mdr F=m+1$, as we have claimed. To show now that $C''$ is free with the given exponents, it is enough to apply Theorem \ref{thmMich}.
\begin{rk}
\label{rkT4}
 If one likes to use Theorem \ref{thmCTC} as above to prove that the curve $C'$ is free, one needs to compute the total Tjurina number $\tau(C')$. This in turn is complicated, since the singularities of $C'$ at the points $p_1=(1:0:0)$ and $p_3=(0:0:1)$ are no longer quasi homogeneous, and hence
$\tau(C'',p_j) <\mu(C'',p_j)$, for $j=1$ and $j=3$. In fact, our Theorem \ref{thm4} combined with Theorem \ref{thmCTC} implies that
$$\tau(C'',p_1) =\tau(C'',p_3) =5m^2-2m,$$
perhaps a result that would not be easy to prove otherwise.
\end{rk}

\subsection{The proof of Theorem \ref{thm5}} 
By its very construction, it is clear that $p$ is a modular point for $C'$.
Let $e=\deg C$ and $m=|\A|=\mult_p(C')$. Hence $d= \deg(C')=e+m$.
Then the fibration
$$\pi_p : M(C') \to B$$
induced by the central projection with center $p$ has as a fiber $F$ the projective line $\PP^1$ minus $e+1$ points, and as a base $B$ the projective line $\PP^1$ minus $m$ points. It follows that the Euler number $E(M(C'))$ of the complement $M(C')$ is given by
$$E(M(C'))=E(F)E(B)=(1-e)(2-m).$$
On the other hand, we know that
$$E(M(C'))=E(\PP^2)-E(C')=3-(\mu(C')-d(d-3)),$$
where $\mu(C)$ is the total Milnor number of $C$.
The above two equations give us
$$\mu(C')=(e+m)^2-em-2m-e+1.$$
Since all the singularities of $C'$ are supposed to be quasi homogeneous, we get the following equality for the total Tjurina number $\tau(C')$ of $C'$
$$\tau(C')=\mu(C')=(e+m)^2-em-2m-e+1.$$
To show  that $C':F=0$ is free we apply \cite[Theorem 1.10]{Mich}.
It follows that $r= \mdr F$ satisfies one of the following properties.
\begin{enumerate}
\item[a)] $r=e$. Then 
$$\tau(e+m,e)_{max}=(e+m-1)^2-e(m-1)=\tau(C'),$$
which implies that $e<m$ and $C'$ is free with exponents $(e,m-1)$ using Theorem \ref{thmCTC}.
\item[b)] $r=m-1$. Then
$$\tau(e+m,m-1)_{max}=(e+m-1)^2-(m-1)e=\tau(C'),$$
which implies that $m \leq e$ and $C'$ is free with exponents $(m-1,e)$ again by Theorem \ref{thmCTC}.
\item[c)] $m\leq r < e.$ This case is impossible, since it implies
$$\tau(C') \leq \tau(e+m,r)_{max}<\tau(e+m,m-1)_{max}=\tau(C').$$
The first inequality follows from Theorem \ref{thmCTC}, and the second one by the fact that the function $t \mapsto \tau(e+m,t)_{max}$ is decreasing for $2t<e+m$.
\end{enumerate}
To prove the last claim in Theorem \ref{thm5}, first notice that the possible non quasi homogeneous singularities of $C'$ may occur only at the intersection $s=L \cap C$, where $s$ is a singular point of $C$ and $L$ is a line in $\A$. If $\mult_sC=2$ and $L$ is not in the corresponding tangent cone $TC_s(C)$, as we have assumed, then
$\mult_sC'=3$ and the $3$-jet $j^3g$ of a local equation $(C',s):g=0$
is a binary cubic form with at least $2$ distinct factors. It follows from the classification of singularities, see for instance \cite{RCS}, that such a singularity has type $D_k$, for some $k \geq 4$, and in particular it is quasi homogeneous.
If $\mu(C,s)=1$ and  $L$ is  in the corresponding tangent cone $TC_s(C)$, then $(C,s)$ is a node $A_1$, and the same argument as above works, namely $(C',s)$ is a $D_k$ singularity. Finally, when
$\mu(C,s)=2$ and  $L$ is  in the corresponding tangent cone $TC_s(C)$, then $(C,s)$ is a cusp $A_2$, and the new singularity
$(C',s)$ is easily seen to be of type $E_7$, hence again quasi homogeneous.

\subsection{The proof of Proposition \ref{prop5}} 

 The point
 $p=(1:0:0)$ is a modular point in this case since any line $L_p$ through $p$,
 not an irreducible component for $C''$, is given by $z=ty$ with $t \ne 0$ and $t^m+1\ne 0$. The intersection $L_p \cap C''$ is described by the equation 
 $$ty^2(x^my^m+t^my^{2m}+t^mx^my^m)(y^m+t^my^m)=t(t^m+1)y^{2m+2}((1+t^m)x^m+t^my^m)=0.$$
 The solution $y=0$ corresponds to the point $p$, which has multiplicity $2m+2$ on $C''$, and there are $m=\deg F'' - (2m+2)$ other intersection points coming from the solutions of $(1+t^m)x^m + t^my^m=0$.
 \begin{ex}
\label{exKP}
Here we show first that the free curves $C'$ and $C''$ coming from Theorem \ref{thm3} are not supersolvable.  For the curve $C'$, it is clear that the only candidates for modular points are the point $p=(0:0:1)$ and the 2 other points obtained from $p$ by permutation of coordinates.
Indeed, a modular point has to contain all the tangents to the Fermat curve issued from it. Now the point $p$ is not a modular point for $C'$, since the line $L_p:y-x=0$ is not an irreducible  component of $C'$ and it satisfies
$$2d+2=|L_p \cap C'| < |L' \cap C'|=\deg C' - \mult_p(C')+1=4d+3-(d+2)+1=3d+2.$$
The same line $L_p:y-x=0$ shows that the point $p$ is not a modular point for $C''$ either. The point $p'=(1:0:0)$ is also not a modular point, as the choice of the line $L_{p'} : z=0$ shows. The point $p''=(0:1:0)$ has the same property, as our curve is invariant under the coordinate change
$x \mapsto y$ and $y \mapsto x$.
To show that the free curve $C'$  coming from Theorem \ref{thm4} is not supersolvable, we use the same approach as above, the lines $L_p$ to use in this case are given by $x=0$, $y=0$ or $z=0$, respectively.
 \end{ex}
\section{Examples: the case of smooth quartic curves} 

In this section we discuss examples of smooth quartic curves having the maximal possible number of flexes of high order.
Let us recall that by Theorem \ref{thmA} the maximal possible number of flex points of order $2$ for smooth quartics is $12$. It is natural to wonder whether there exists a complete classification of smooth quartics which have exactly $12$ flexes of order $2$. In order to do so, we discuss interesting properties of the following pencil of quartics which was studied by Ciani in the 19th century.

Let us define
\begin{equation}
\label{quarticpencil}
C_{\lambda} \, : \, x^{4} + y^{4} + z^{4} + \lambda\cdot (y^{2}z^{2} + z^{2}x^{2} + x^{2}y^{2})=0.
\end{equation}
It is easy to observe that each curve in the pencil is invariant under the natural action of an octahedral group of collineations.

There are some values of $\lambda$ which lead to special members of the pencil, namely
\begin{itemize}
\item $\lambda = 0$ gives us the Fermat quartic curve, or Dyck's curve, which has a group of $96$ collineations;
\item if $\lambda$ is a root of $\lambda^{2}+3\lambda+18=0$, then we get the Klein quartic curve having a group of $168$ collineations.
\end{itemize}
In the case of the Fermat quartic, by a discussion presented above, we know that it has exactly $12$ flexes of order $2$, so the maximal possible number in the class of smooth quartics. In the case of the Klein quartic curve, we know that this curve has only flex points of order $1$, so exactly $24$ flexes. Now we pass to another interesting element in the pencil of quartics by taking $\lambda=3$. The resulting quartic $C_{3}$ is smooth and it has the group of collineations of order $24$. It was verified directly by Edge in \cite{Edge} that the curve $C_{3}$ admits exactly $12$ flexes of order $2$ and he provided both the coordinates of these points and the equations of the associated tangent lines. Now we recover Edge's calculations. Looking precisely on the Hessian $H$ of $C_{3}$, which is
$$H=2x^6 + x^{4}(3z^2+3y^{2})+x^{2}(8y^{2}z^{2}+3z^{4}+3y^{4})+2z^{6}+2y^{6}+3z^{4}y^{2}+3z^{2}y^{4},$$
one can show that flexes of order $2$ are just the intersection points of the curve $C_{3}$ with the $6$ lines given by the linear factors of
\begin{equation}
\label{eqKK}
F=(x^{2}+y^{2})(y^{2}+z^{2})(z^{2}+x^{2}).
\end{equation}
The flex points have the following coordinates:
\[ \begin{tabu} {ll}
P_{1} : (i:1:-1), & P_{2} : (-i:1:-1), \\
P_{3} : (-1:i:1), & P_{4} : (-1:-i:1), \\
P_{5} : (1:-1:i), & P_{6} : (1:-1:-i), \\
P_{7} : (-i:1:1), & P_{8} : (i:1:1), \\
P_{9} : (1:-i:1), & P_{10}: (1:i:1), \\
P_{11}: (1:1:-i), & P_{12}: (1:1:i). 
\end{tabu}\]
Observe that these $12$ flexes of order $2$ are uniformly distributed, four on each of the lines defined by $F$.  

Up to now we described exactly two smooth quartics having the maximal possible number of flexes of order $2$. However, as it turns out by a result due to Kuribayashi and Komiya \cite{KK}, these are the only smooth plane quartic curves having $12$ flexes of order $2$, and this is rather surprising.

\begin{rk}
\label{rkC3}
We have seen in Theorem \ref{thm3} that if we add to the Fermat quartic its $12$ inflectional tangents of order $2$ and the triangle $\Delta:xyz=0$ determined by the inflection points, then we get a free curve of degree 
$19$. If we try to apply the same construction to the quartic curve $C_3$, the resulting curves are far from being free. One explanation for this fact may be the following. The union $\A_F$ of the $12$ inflectional tangents 
of the Fermat quartic is a line arrangement having $3$ points of multiplicity $4$. On the other hand, the union $\A_3$ of the $12$ inflectional tangents 
of the quartic $C_3$ is a line arrangement having only double points, and hence the total Tjurina number $\tau(\A_3)$ is much smaller than $\tau(\A_F)$. If we add the triangle $\Delta$ to $\A_F$, we get a line arrangement having 3 points of multiplicity $6$. On the other hand, if we add to $\A_3$ the 6 lines determined by \eqref{eqKK}, we get a line arrangement having only points of multiplicity $2$ and $3$, and hence having small total Tjurina number compared with respect to its degree.
\end{rk}

\section{Examples: the case of singular curves} 

\begin{ex}
\label{ex1}
Any nodal cubic is projectively equivalent to the cubic
$$C:F=xyz+x^3+y^3=0.$$
The corresponding Hessian is $H= -2(3(x^3+y^3)-xyz)$. Hence the intersection $C \cap H_C$ consists of the following 4 points:
$$p_1=(0:0:1) \text{ and } p_{j}=(1:j:0),$$
where $j^3+1=0$. The point $p_1$ is a simple node, and the points $p_j$ give rise to $3$ inflection points of order $1$. This is reflected in the equality
$$(C,H_C)_{p_1} +\sum_j (C,H_C)_{p_j}=6+1+1+1=9,$$
recall Example \ref{exA1}. 
The $3$ inflectional tangents $L_k$ for $k=1,2,3$ are given by the equations $L_k: 3x+3j_k^3y+j_kz=0$, where $j_k$ are the 3 roots of the equation $j^3+1=0$. It follows that these $3$  inflectional tangents $L_k$ are not concurrent, so their addition to $C$ will not give free curves as in Remark \ref{rkC3} above.
On the other hand, if we add to $C$ the tangent cone at the singular point, we obtain the curve
$$C': F'=xy(xyz+x^3+y^3)=0,$$
which is free with exponents $(2,2)$ as a direct computation with \verb}SINGULAR} shows. Moreover, this curve $C'$ is supersolvable, since clearly $p_1$ is a modular point for $C'$.

Any cuspidal cubic is projectively equivalent to the cubic
$$C:F=x^2z+y^3=0.$$
The corresponding Hessian is $H= -24x^2y$. Hence the intersection $C \cap H_C$ consists of the following 2 points
$$p_1=(0:0:1) \text{ and } p_2=(1:0:0).$$
The point $p_1$ is a cusp $A_2$, and the point $p_2$ is an inflection point of order $1$. This is reflected in the equality
$$(C,H_C)_{p_1} + (C,H_C)_{p_2}=8+1=9,$$
recall Example \ref{exA3}. This is a special case of Theorem \ref{thm2}, and gives rise to two free curves by adding one or two lines, as explained there.
\end{ex}

\begin{ex}
\label{ex2}
In this example we consider some plane quartic curves. 

Consider  the quartic $C: F=(x^3+y^3)z+x^4+y^4=0$, which has a $D_4$-singularity at $p_1=(0:0:1)$. The corresponding Hessian is
$$H=-54(xyz(x^3+y^3)+2x^2y^2(x^2+y^2)).$$
The point $p_1$ is an ordinary simple singularity of multiplicity $k=3$, and hence $(C,H_C)_{p_1}=18$ by
Theorem \ref{thm1}. There are in addition $6$ inflection points of order $1$, with are the points $(1:0:-1)$, $(0:1:-1)$ and the $4$ points $(u:v:w)$, where $(u:v)$ is coming from the $4$ solutions of the equation
$$u^4+v^4-2uv(u^2+v^2)=0$$
in $\PP^1$ and $w=-(u^4+v^4)/(u^3+v^3)$. If we add the tangent cone of the singular point, namely the lines $x^3+y^3=0$, we get a free curve 
$$C':F'=(x^3+y^3)F=0,$$
of degree $7$ and exponents $(3,3)$.  Moreover, this curve $C'$ is supersolvable, since clearly $p_1$ is a modular point for $C'$.

Next, consider the quartic
\begin{equation}
\label{eq3A1}
C:F=x^2y^2+y^2z^2+x^2z^2=0,
\end{equation}
which has $3$ nodes. It is easy to see that all of them have type $(3,3)$, and hence $C$ has no inflection points by Example \ref{exA1}. The corresponding Hessian is
$$H=-24(x^4y^2+x^2y^4+y^4z^2+y^2z^4+x^4z^2+x^2z^4-6x^2y^2z^2).$$
If we add to $C$ one tangent line at each of  the $3$ nodes, namely the lines
$$(x+iy)(y+iz)(z+ix)=0,$$
we get a free curve of degree $7$ with exponents $(3,3)$. All the singularities of this curve are simple, but this curve is not maximizing, recall our discussion in Section 3 on these curves.
Finally the quartic
$$C:F =x^2y^2+y^2z^2+x^2z^2-2xyz(x+y+z)=0,$$
which has $3$ cusps $A_2$. Hence $C$ has no inflection points by Example \ref{exA3}. The corresponding Hessian is
$$H=144(x^3y^3+y^3z^3+x^3z^3-x^3(y^2z+yz^2)-y^3(x^2z+xz^2)-z^3(x^2y+xy^2).$$
Let $L_1$, $L_2$ and $L_3$ be the 3 lines which are the reduced tangent cones corresponding to the $3$ cusps, which are given up to an order by the equations $x-y=0$, $y-z=0$ and $z-x=0$. Then the curves
$$C_1=C \cup L_1, \ C_2=C_1 \cup L_2 \text{ and } C_3=C_2 \cup L_3$$ are all free, with exponents respectively
$$ (2,2), \ (2,3) \text{ and } (2,4).$$
We get in this way a free curve $C_1$ of degree $5$ and maximizing curves $C_2$ and $C_3$, of degree $6$ and $7$, respectively, as already pointed out in \cite{MaxCurv}. It is interesting to note that the curve $C_3$ is supersolvable, and the point $p=(1:1:1)$ is a modular point for it. Indeed, the lines joining $p$ to the singularities of $C$ are already in $C_3$. It remains to show that any line $L_p$ through $p$, different from $L_1,L_2,L_3$ meets $C$ in exactly 4 points, that is $L_p$ is not a tangent line to $C$. If $q=(u:v:w) \in C$ is a smooth point such that the tangent line $T_qC$ passes through $p$, then we have
$$F_x(q)+F_y(q)+F_z(q)=0.$$
A direct computation shows that
$$F_x(q)+F_y(q)+F_z(q)=-12uvw$$
and hence at least one of the coordinates of $q$ vanishes.
But then $q \in C$ implies that 2 coordinates vanish, and therefore $q$ is a singularity of $C$, a contradiction. Note that the curve $C_3$ is an example of curve satisfying both the first assumption in Theorem \ref{thm5},
since all of its singularities are quasi homogeneous, and the second assumption, even if $p$ belongs to the tangent cones $TC_s(C)$ of the three cusps, as they have Milnor numbers equal to $2$.
\end{ex}

\begin{ex}
\label{ex3}
In this example we consider the curve
$$C: F=x^my^m+y^mz^m+x^mz^m=0$$
of degree $d=2m \geq 4$. This curve has 3 ordinary $m$ multiple points
$$p_1=(1:0:0), \ p_2=(0:1:0) \text{ and } p_3=(0:0:1)$$
which have type $\underbrace{(m+1, \ldots, m+1)}_{m \text{ times }}$. These singularities are easily seen to be quasi homogeneous and hence
$$\mu(C,p_j)=\tau(C,p_j)=(m-1)^2,$$
for $j \in \{1,2,3\}$. It follows by Theorem \ref{thmA} that
$$(C,H_C)_{p_j}=3(m-1)^2+m-3+m(m+1)=4m(m-1).$$
The equality \eqref{eq2.1} implies that $C$ has no inflection points. The reader can check that this curve $C$ is not even nearly free, e.g. for $m=3$. On the other hand, we show now that the curve
$$C':F'=xyzF=0$$
is free with exponents $(m+1,m+1)$. In order to show this we note first that the only singularities of $C'$ are again the points $p_j$ for $j \in \{1,2,3\}$,
which are ordinary quasi homogeneous singularities of multiplicity $(m+2)$. For the last claim one can use \cite[Exercise (7.33)]{RCS}.
It follows that
$$\tau(C')=3(m+1)^2.$$
The equality
$$x(y^m-z^m)F_x-y(y^m+z^m)F_y+z(y^m+z^m)F_z=0$$
shows that $\mdr (F)=m+1$. This implies that one has $\mdr F' \geq \mdr F=m+1$ and
$$\tau(2m+3,m+1)_{max}=(2m+2)^2-(m+1)^2=\tau(C').$$
This equality implies our claim by Theorem \ref{thmCTC}.
\end{ex}

\begin{ex}
\label{ex4}
In this example we construct a {\it new example of maximizing curve in degree 8 and the first and unique known example of maximizing curve in degree} $9$.
 Such maximizing curves of odd degree seem to be quite exceptional, for instance we know no example in degrees $\geq 11$.
The curve
$$C:F=(x^2+y^2+z^2)^3-27x^2y^2z^2=0$$
is the dual of the quartic with 3 nodes in \eqref{eq3A1}.
This sextic curves has  six cusps  $A_2$ at the points $(0:1:\pm i)$, $(1:0: \pm i )$ and $ (1:\pm i: 0)$, where $i^2=-1$ and four nodes $A_1$ at the points $(1:\pm 1: \pm 1)$. This curve is far from being free,
since the corresponding module ${\rm D}_{0}(F)$ has 4 generators, as one can check using, for instance, \verb}Singular}. 
If we add the line $L_1:x=0$, the curve $C_1=C \cup L_1$ is nearly free with exponents $(3,4)$. If we add one more line, namely $L_2:y=0$, the resulting curve $C_2=C_1 \cup L_2$ is free with exponents $(3,4)$, and this gives a {\it new example of maximizing curve in degree} 8.
Finally, if we add the third line $L_3:z=0$, the resulting curve
$$C_3=C_2 \cup L_3: F'=xyz\left( (x^2+y^2+z^2)^3-27x^2y^2z^2 \right)=0$$
has six singularities $E_7$ at the points $(0:1:\pm i)$, $(1:0: \pm i )$ and $ (1:\pm i: 0)$, and seven nodes $A_1$ at the points
$(1:\pm 1: \pm 1)$, $(1:0:0)$, $(0:1:0)$ and $(0:0:1)$.
It follows that
$$\tau(C_3)=6\cdot 7+7 \cdot 1=49,$$
and this equality implies our claim by Definition \ref{defMC}. Theorem \ref{thmMC} tells us that $C_3$ is a free curve with exponents $(3,5)$ obtained from $C$ by adding the three lines $L_1$, $L_2$ and $L_3$.
\end{ex}

\section{On a theorem by Schenck, Terao and Yoshinaga and conic-line free arrangements} 

We start with a remark concerning the paper \cite{STY}. Let $C_1:F_1=0$ and $C_2:F_2=0$ be two reduced curves in $\PP^2$, without common irreducible components. We denote $d_j=\deg F_j$ and $r_j=mdr(F_j)$ for $j=1,2$. Consider now the union of the two curves
$C:F=F_1F_2=0$, and let $d=d_1+d_2=\deg F$ and $r= {\rm mdr}(F)$. Using the main result in \cite{STY}, namely  \cite[Theorem 1.6 and Remark 1.8]{STY}, one can obtain relations among the 3 integers $r_1,r_2$ and $r$.
The hypothesis in  \cite[Theorem 1.6 and Remark 1.8]{STY} are the following: all the singularities of $C_1$ and $C$ are quasi homogeneous, and $C_2$ is a smooth curve. However, the quasi homogeneity of the singularities of $C_1$ and $C$ is only used in
 \cite[Proposition 2.5]{STY} to compute the difference
 $$\tau(C)-\tau(C_1)$$
 using  \cite[Lemma 2.4]{STY}, which is the equality \eqref{eqMilnor} above. In this difference, the contribution of the singularities of $C_1$ not on $C_2$ cancels, and hence we need to control only the change in the Tjurina numbers at the points of the intersection $C_1 \cap C_2$, when we add the curve $C_2$. Hence the hypothesis in  \cite[Theorem 1.6 and Remark 1.8]{STY} can be relaxed to the following: all the singularities of $C_1$ and $C$ {\it situated on $C_2$} are quasi homogeneous, and $C_2$ is a smooth curve. 
 In view of this remark, we have the following stronger version of \cite[Corollary 6.4]{mdr}.

\begin{cor}
\label{corDIS}
 With the above notation, assume that $C_2$ is a smooth curve, and that all singularities of the curves $C_1$ and $C$ situated on $C_2$ are quasi homogeneous.  Let $R$ be the reduced scheme of $C_1 \cap C_2$. If
$$|R|>(r_1+1)d_2,$$
then $r=r_1+d_2$.
\end{cor}

We use this stronger result to construct a family of free conic-line arrangements as follows. Let $C_0$ be the union of $m$ smooth conics belonging to a hyperosculating pencil of conics, that is a pencil of conics with one base point. Let $C_1=C_0 \cup L_0$, where $L_0$ is the common tangent to all the conics in $C_0$.
An explicit equation for a special case of such curves $C_0  :F_0=0$ and $C_1 : F_1=0$ can be found in 
\cite[Equation (1.8)]{DK}, namely
\begin{equation}
\label{eqDIS1}
F_0=x^{2m}+(xz+y^2)^m \text{ and } F_1 = xF_0 = x(x^{2m}+(xz+y^2)^m).
\end{equation}
These curves have been considered by several authors, see for instance \cite[Section 7.5, p. 179]{CTC} and \cite{Pl,Shin}, in relation with the maximal Milnor number a singular point on a plane curve of degree $d$ might have.
These authors showed in particular that
$$\mu(C_j)=(d_j-1)^2- \bigg\lfloor \frac{d_j}{2} \bigg\rfloor,$$
where $d_j=\deg C_j$ for $j=0,1$. On the other hand, it is known that these curves are free with $r_j= {\rm mdr}(F_j)=1$, for $j=0,1$, see \cite[Theorem A]{Be}, where $C_0$ is denoted by $\CC_1$ and $C_1$ is denoted by $\CC\LL_1$. See also \cite{DK} for the particular curves in \eqref{eqDIS1}. It follows that
$$\tau(C_j)=(d_j-1)^2-d_j+2 <\mu(C_j)$$
as soon as $d_j \geq 5$. Hence the singularity of the curve $C_j$ at the point $p=C_0 \cap L_0$, the base point of the pencil, is not quasi homogeneous when $d_j\geq 5$. Let $q \in L_0$ be a point distinct from $p$. Let $Q_1,\ldots, Q_m$ be the smooth conics in $C_0$ and let $L_j$ be the tangent from $q$ to the conic $Q_j$, distinct from the tangent $L_0$. With this notation we have the following.

\begin{prop}
\label{propDIS1}
For any $m \geq 2$ and $j=2,\ldots, m+1$, the conic-line arrangement
$$C_j=C_1 \cup L_1 \cup \ldots \cup L_{j-1}: F_j=0$$
is a free curve of degree $d_{j}=2m+j$ with $r_{j}=j$.
\end{prop}
\proof
We prove this claim by induction on $j$. The case $j=1$ is proved in \cite{Be}, as we said above. Assume that the claim holds for some curve $C_k$, with $1\leq k<m+1$, and let's prove it for $C_{k+1}$. The fact that $d_{k+1}=2m+k+1$ is obvious. Note that the line $L_{k}$ is tangent to the conic $Q_k$ and meets any other conic $Q_j$ for $j \ne k$ in 2 points. To see this, one may use Remark \ref{rkDIS0} below. Hence
$$|R_{k}|=|C_k \cap L_{k}|=2(m-1)+1+1=2m>k+1=r_k+1,$$
since $k<m+1 \leq m+1+(m-2)=2m-1$. Using Corollary \ref{corDIS} we get $r_{k+1}=r_k+1=k+1$. Note that
$$2r_{k+1} =2(k+1) \leq d_{k+1}-1=2m+k$$
since $k \leq 2m-2$ as we have seen above. Hence, using Theorem \ref{thmCTC}, the curve $C_{k+1}$ is free if and only if
$$\tau(C_{k+1})=\tau(d_{k+1},r_{k+1})_{max}.$$
In other words, we have to check that the addition of the line $L_k$ to
$C_k$ increases its total Tjurina number by
$$\Delta_k=\tau(d_{k+1},r_{k+1})_{max}-\tau(d_{k},r_{k})_{max}=
d_k+r_k=2m+2k.$$
At the point $q$, the change in Tjurina number is
$$k^2-(k-1)^2=2k-1,$$
since we pass essentially from a local equation $u^k+v^k=0$ to
$u^{k+1}+v^{k+1}=0$. Except the point $q$, there is a singularity of type $A_3$ corresponding to the tangent point, and $2(m-1)$ nodes $A_1$ on the line $L_k$.
Hence, we get
$$2k-1+3+2(m-1)=\Delta_k,$$
which completes our proof.
\endproof

\begin{rk}
\label{rkDIS0}
Note that a line $L_k$, which is tangent to the conic $Q_k$, cannot be also tangent to another conic $Q_j$ with $j \ne k$. To show this, it is enough to consider the curve $C_1$ when $m=2$. If there is a point $q \in L_0$ such that the equality $L_1=L_2$ holds, then the curve
$C_2$ has degree $d_2=6$, $r_2=2$ by Corollary \ref{corDIS}, and the contribution of the line $L_1$ to $\tau(C_2)$ would be
$$ 3+3 +1 >\Delta_1=5+1.$$
This is a contradiction with Theorem \ref{thmCTC} since our calculation above would yield
$$\tau(C_{2}) > \tau(6,2)_{max}.$$
\end{rk}
\begin{rk}
\label{rkDIS1}
Note that the curve $C_{m+1}$ is supersolvable since $q$ is a modular point. Since the singularity at $p$ is not quasi homogeneous for $m\geq 2$, we get a new countable family of supersolvable curves, supporting our Conjecture \ref{conj1}, and not covered by Theorem \ref{thm5}. An explicit equation for such a curve is the following, obtained from \eqref{eqDIS1}:
$$C_{m+1}: F_{m+1}=x(x^m+z^m)(x^{2m}+(xz+y^2)^m)=0,$$
where $p=(0:0:1)$ and $q=(0:1:0)$.
\end{rk}

\begin{rk}
\label{rkDIS2}
In \cite[Theorem A]{Be}, the authors classify all the conic-line arrangements $C:F=0$ which are free with ${\rm mdr}(F)=1$.
Besides the union of concurrent lines and the curves $C_0$ and $C_1$ discussed above, the other 4 types of such curves, denoted by $\CC\LL_j$ with $j=2,3,4,5$, are obtained using a bitangent pencil of conics. It is easy to see, using for instance \cite[Equation 5.1]{mdr} with $x^2-y^2$ replaced by $xy$, that in these cases all the singularities of these curves are quasi homogeneous. One can obtain new families of free conic-line arrangements from them by adding concurrent tangent lines, but, in particular, the resulting supersolvable curves are covered by our Theorem \ref{thm5}. As an example, consider
the  curve $ \mathcal{CL}_5$ studied in \cite{Be}. This curve, which we denote by  $D_0$, is  the union of $m$ smooth conics belonging to a bitangent pencil of conics, the tangent lines $T_0$ and resp. $T_0'$ in the tangency points $p$ and resp. $p'$ of the pencil, and the line $L_0$ joining the two tangency points, see
\cite[Theorem A]{Be}. Let $q$ be a point on $L_0$, distinct from $p$ and $p'$. Then add one by one tangents to the smooth conics in $D_0$ passing through $q$. For each conic we can add one or two tangents. If we denote by $D_k$, for $1 \leq k \leq 2m$, one of the several conic-line arrangements obtained from $D_0$ by adding $k$ such tangent lines, it it easy to show exactly as in the proof of Proposition \ref{propDIS1} that
$D_k:F_k=0$ is a free conic-line arrangement with
$d_k=\deg F_k= 2m+3+k$ and $r_k=k+1$.  Note that in this case the difference $\delta_k$ given by
$$\delta_k=d_k-1-2r_k=2m-k \geq 0$$
can be as small as we like, can even vanish when $k=2m$, that is for the  curve corresponding to $k=2m$ in this family. On the other hand, for the curves $C_k$ in Proposition \ref{propDIS1} one has
$$\delta_k=d_k-1-2r_k=2m+k-1-2k =2m-k-1\geq m-2 >0$$
when $m>2$. 

To get a supersolvable curve from $D_{2m}$, it is enough to add the line $L_0'$ joining the point $q$ with the intersection point $q'=T_0 \cap T_0'$. Note that 
$$|D_{2m}\cap L_0'|= 2m+2=r_{2m}+1,$$
and hence Corollary \ref{corDIS} cannot be applied. On the other hand, we know that $D_{2m+1}=D_{2m}\cup L_0': F_{2m+1}=0$ is a free curve by Theorem \ref{thm5}, and hence
$$2m+1=r_{2m}\leq r_{2m+1} \leq \frac{\deg F_{2m+1}-1}{2}=\frac{4m+3}{2}.
$$It follows that one has 
$$2m+1=r_{2m}= r_{2m+1} .$$
\end{rk}

\end{document}